\documentclass[onecolumn,12pt]{IEEEtran}
\usepackage{graphicx} % Required for inserting images

\usepackage{amsmath}
\usepackage{amssymb}
\usepackage{hyperref}

\usepackage{tikz}
\usetikzlibrary{matrix,positioning,calc}

\usepackage{setspace}
\usepackage{comment}
\usepackage{url}
\usepackage{hyperref}

\usepackage{mleftright}

\newcommand{\such}{\, | \,}

\newcommand{\tr}{\operatorname{tr}}

\newcommand{\diag}{\operatorname{diag}}

\newcommand{\be}{\begin{equation}}
\newcommand{\ee}{\end{equation}}
\newcommand{\R}{\mathbb R}
\newcommand{\C}{\mathbb C}

\DeclareMathOperator{\vol}{volume}
\DeclareMathOperator{\vect}{vec}
\DeclareMathOperator{\vech}{vech}
\DeclareMathOperator{\ind}{index}

\DeclareMathOperator{\image}{Im}

\newtheorem{Theorem}{Theorem}
\newtheorem{Lemma}{Lemma}

\newtheorem{Corollary}{Corollary}
\newtheorem{Definition}{Definition}
\newtheorem{Remark}{Remark}
\newtheorem{Example}{Example}

\newcommand{\updt}[1]{\textcolor{blue}{#1}}

\makeatletter
\newcommand{\ostar}{\mathbin{\mathpalette\make@circled\pi}}
\newcommand{\make@circled}[2]{%
  \ooalign{$\m@th#1\smallbigcirc{#1}$\cr\hidewidth$\m@th#1#2$\hidewidth\cr}%
}
\newcommand{\smallbigcirc}[1]{%
  \vcenter{\hbox{\scalebox{0.77778}{$\m@th#1\bigcirc$}}}%
}
\makeatother

\title{Multiplicative and additive compounds via  Kronecker products and Kronecker sums\thanks{This research was partly supported by research grant 407/19 from the~Israel Science Foundation.  } }
\author{Ron Ofir and Michael Margaliot\thanks{RO is with the Andrew and Erna Viterbi Faculty of Electrical and Computer Engineering, Technion---Israel Institute of Technology, Haifa 3200003, Israel. MM (corresponding author) is with the School of Electrical Engineering, Tel Aviv University, Israel 69978 (e-mail: michaelm@tauex.tau.ac.il)
}}

\begin{document}

\maketitle
\doublespacing
\begin{center}
    {31 December  2023}
\end{center} 

\begin{abstract}
%%%%%%%%%%%%%%%%%%%%

Compound matrices play an important role in many fields of mathematics and have recently found new applications in systems and control theory. However, the explicit formulas for these compounds are non-trivial and not always easy to use. 
Here, we derive new formulas for the multiplicative and additive compounds of a matrix using Kronecker products and sums. This provides a new approach to matrix compounds based on the well-known and powerful theory of Kronecker products and sums. 
We demonstrate several applications of these new formulas, including deriving a new expression for the additive compound of the product of two matrices. 
%%%%%%%%%%%%%%%%%%%%%%
%%%%%%%%%%%%%%%%%%%%%%
\end{abstract}


\begin{center}
\end{center}

 \section{Introduction}
%%%%%%%%%%%%%%%%%%%%%%%%%

Compound matrices play an important role in linear algebra, geometry, graph theory, and dynamical systems (see, e.g.,~\cite{fiedler_book,Gantmacher_vol1,muldowney1990compound,BAPAT_minor_laplacian,gk_book,Wang2add,Fomin1999TotalPT}). 
Recently, they have found interesting  applications in systems and control theory (see the tutorial paper~\cite{comp_long_tutorial}). One reason for this is that $k$-compound matrices allow to track the volume of a parallelotope spanned by~$k$ vectors under the evolution of an $n$-dimensional  time-varying linear dynamics. Following the seminal work of Muldowney~\cite{muldowney1990compound}, this allows to generalize the theory of contractive systems~\cite{LOHMILLER1998683,sontag_cotraction_tutorial} to $k$-contractive systems, with~$k$ an integer,~\cite{kordercont} and~$\alpha$-contractive systems, with~$\alpha$ real~\cite{wu2020generalization}. These generalizations have been used to analyze the global behaviour of nonlinear dynamical systems in several fields of science including recurrent neural networks~\cite{Ron_Lurie}, mathematical epidemiology~\cite{kordercont}, chaotic systems~\cite{wu2020generalization}, and systems biology~\cite{ron_general_lurie}.

Compound matrices also play an important role in totally positive differential systems~\cite{margaliot2019revisiting}, and~$k$-cooperative systems~\cite{Eyal_k_posi}.
For more on the applications of compound matrices in systems and control theory, see the tutorial paper~\cite{comp_long_tutorial}, and the recent papers~\cite{Angeli2014hopf,Angeli2023smallgain,Angeli2021nonosc,grussler2021internally,k_comp_DAE,LMI_Kcont}.

However, the known explicit expressions for compound matrices are non-trivial and this may hamper important applications.
For example, the expressions for compounds of block matrices destroy the block form. If the block form represents a system composed of the interconnection of several sub-systems then this makes it difficult to understand how each sub-system affects the compounds.

%%%%%%%%%%%%%%%%%%

Attempts to better understand  the structure of compound matrices date back at least to the work of Aikten~\cite{Aikten_normal_form}, who used a combinatorial  approach to relate the Jordan form of a matrix to the Jordan form of its compounds (see also~\cite{aitken_book}). More recently, the fact that the sufficient condition for~$k$-contraction of a nonlinear system  is based on the~$k$-additive compound of the Jacobian of the vector field has led to a renewed interest in the compounds of block matrices~\cite{ofir2021sufficient, Angeli2023smallgain,angeli_bounding_lyap}.

Here, we derive new formulas for the $k$-multiplicative and $k$-additive compounds of matrices in terms of Kronecker products and sums, see Theorems~\ref{thm:multi_compi}
and~\ref{thm:kron_addi_compound}    below.   Roughly speaking, this allows to 
re-derive  the   entire theory of   multiplicative and additive compounds, of any order~$k$,
using the well-known and powerful theory of Kronecker products and sums. We demonstrate 
this using several examples. 
  We also derive  
several new applications of our new formulas. In particular, we derive a new formula for the~$k$-additive compound of the product of two matrices. 
Some of the theoretical
results in this  paper can be derived using the framework  of multilinear algebra and wedge products, but we  opt  an approach that is based on matrix theory.

We use standard notation. Vectors [matrices] are denoted by small [capital] letters. We focus on   real matrices, as this is the relevant case in applications. For a vector~$x\in\R^n$, $|x|_2:=(\sum_{i=1}^n x_i^2)^{1/2}$ is the~$L_2$ norm of~$x$.  
 The transpose of a matrix~$A$ is~$A ^T$. For a square matrix~$A$, $\det(A)$ [$\tr(A)$] denotes the determinant [trace] of~$A$.  The image of a matrix~$A\in\R^{n\times m}$  is the subspace 
\[
\image A :  =\{Ax \such  x\in\R^m\}.
\]
We use~$I_{n}$ to denote the~$n\times n$ identity matrix. 
For a matrix~$A\in\R^{m\times n}$, the vectorization of~$A$, denoted~$\vect(A)$, is the column vector  with $mn$ entries  obtained by stacking the entries of~$A$ column after column. 
For example,
for~$A\in\R^{3\times 2}$, we have
\[
\vect\mleft(\begin{bmatrix}a_{11}& a_{12} \\
a_{21}& a_{22} \\ a_{31} & a_{32} 
\end{bmatrix} \mright) =\begin{bmatrix}
    a_{11} &a_{21}& a_{31}& a_{12} &a_{22} & a_{32} 
\end{bmatrix}^T. 
\]
%%%
For an integer~$n\geq 1$ and~$k\in\{1,\dots,n\}$, let~$R(n,k)$~[$Q(n,k)$] denote all the sequences [\emph{increasing} sequences] of~$k$ integers in the range~$\{1,\dots,n\}$, ordered lexicographically. For example, 
\[
    R(3,2) = ( (1,1), (1,2), (1,3), (2,1), (2,2), (2,3), (3,1), (3,2), (3,3) ), 
 \]
and
\[
 Q(3,2) = ( (1,2), (1,3), (2,3) ).
\]
Let~$\ind^R_{n,k}(i_1,\dots,i_k)$ [$\ind^Q_{n,k}(i_1,\dots,i_k)$]
denote the index of the   sequence~$(i_1,\dots,i_k)$ in~$R(n,k)$ [$Q(n,k)$].
For example, 
$
\ind^R_{3,2} (1,3) = 3,
$
and
$
\ind^Q_{3,2} (1,3) = 2. 
$
These functions play an important role in the analysis below, so 
we   detail several  properties of~$\ind^R_{n,k}$ and~$\ind^Q_{n,k} $ in the Appendix.

The remainder of this paper is organized as follows.  For the sake of completeness, the next section reviews 
known 
definitions and results that are used later on. Section~\ref{sec:main} describes the main results. Section~\ref{sec:appli} details several new applications of the theoretical results.  The final section concludes and describes   
possible directions for further research.
A collection of several more technical results are placed in the Appendix.

%%%%%%%%%%%%%%%%%%%%
\section{Preliminaries}\label{sec:prelim}
%%%%%%%%%%%%%%%%%
Our   results  are based on relating
   matrix compounds and Kronecker products/sums. We begin by briefly reviewing these topics. 

%%%%%%%%%%%%%%%%%%%%%%%%%%%
\subsection{Compound  matrices}
%%%%%%%%%%%%%%%%%%%%%%%%%%%%%%%%%%%%%%%
%%%%%%%%%%%%%%%%%%%%%%%
\begin{Definition}\label{def:MC}
%%%%%%%%%%%%%%%%%%%%%%%%%%%
Fix~$A\in\R^{n\times m}$ and~$k\in \{1,\dots,\min\{n,m\}  \} $. The~$k$-multiplicative compound matrix of~$A$, denoted~$A^{(k)}$, is the~$\binom{n}{k}\times\binom{m}{k}$ matrix that includes  all the~$k$-minors of~$A$,
ordered lexicographically.
%%%
\end{Definition}

For example,
for
$
A=\begin{bmatrix}
    a_{11} & a_{12} \\
    a_{21} & a_{22} \\
    a_{31} & a_{32} 
\end{bmatrix}
$
and~$k=2$, we have
$
A^{(2)}=\begin{bmatrix}
   \det  \begin{bmatrix}
   a_{11}& a_{12}\\
   a_{21}&a_{22}
\end{bmatrix} \\
%%%
 \det  \begin{bmatrix}
   a_{11}& a_{12}\\
   a_{31}&a_{32}
\end{bmatrix} \\
%%%
\det  \begin{bmatrix}
   a_{21}& a_{22}\\
   a_{31}&a_{32}
\end{bmatrix} 
%%%%%%%
\end{bmatrix}.
$

Definition~\ref{def:MC}
implies that~$(A^{(k)})^T = (A^T)^{(k)}$,
$A^{(1)}=A$,  and if~$A\in\R^{n\times n}$ then~$A^{(n)}=\det (A)$. 
 %%%%%%
If~$D\in\R^{n\times n}$ is diagonal, that is,~$D=\diag(\lambda_1,\dots,\lambda_n)$ then Definition~\ref{def:MC} implies that~$D^{(k)}$ is also diagonal, with
\[
D^{(k)}=\diag\left(\prod_{i=1}^k\lambda_i, (\prod_{i=1}^{k-1}\lambda_i) \lambda_{k+1},\dots,
\prod_{i=n-k+1}^n\lambda_i\right).
\]
In particular,  
the $k$-multiplicative compound of the~$n\times n$ identity matrix is~$(I_n)^{(k)}=I_r$, with~$r:=\binom{n}{k}$.

More generally, 
the~$k$-multiplicative compound of a square matrix~$A\in\R^{n\times n}$ has a  ``multiplicative''  spectral property. If   the eigenvalues of~$A$ are~$\lambda_1,\dots,\lambda_n\in\C$  then the eigenvalues of~$A^{(k)}$ are  all   the~$\binom{n}{k}$ products:
\[
\lambda_{i_1}\dots\lambda_{i_k}, \quad 1\leq i_1<\dots< i_k\leq n. 
\]

The term ``multiplicative compound'' is justified by the  
Cauchy-Binet theorem that asserts that for any~$A\in\R^{n\times m}$, $B\in\R^{m\times p}$ and~$k\in
\{1,\dots,\min\{n,m,p\}\}$, we have
\be\label{eq:cbin} 
(AB)^{(k)}=A^{(k)} B^{(k)}.
\ee
Note that for~$A,B\in\R^{n\times n}$ and~$k=n$, this reduces to~$\det(AB)=\det(A)\det(B)$. Furthermore, 
if~$A $ is square and invertible then applying~\eqref{eq:cbin} with~$B=A^{-1}$ implies that~$A^{(k)}$ is also invertible, and $(A^{(k)})^{-1} = (A^{-1})^{(k)}$.

\subsubsection{$k$-compounds and~$k$-parallelotopes}
 %%%%%%%%%%%
 
 $k$-compound matrices have a geometric interpretation 
 that plays an important role in  their applications in nonlinear dynamical systems. To explain this,  fix~$k$ vectors~$a^1,\dots,a^k\in\R^n$. The parallelotope generated by these vectors is
\[
\mathcal P(a^1,\dots,a^k):=\left\{\sum_{i=1}^k r_i  a^i\such r_i\in[0,1]\right\}
\]
(see Fig.~\ref{fig:paralellotope}). Define the~$n\times k$ matrix
$
V:=\begin{bmatrix} a^1&\dots& a^k\end{bmatrix} 
$.
Then
\be\label{eq:volupp}
\vol(\mathcal P(a^1,\dots,a^k) )=
|V^{(k)}|_2 
\ee
 (see, e.g.,~\cite{Gantmacher_vol1}). 
Note that~$V^{(k)}$ has dimensions~$\binom{n}{k}\times\binom{k}{k}$, i.e., it is a column vector. 
In particular, if~$a^i=a^j$ for some~$i\not = j$ then~$\begin{bmatrix}
    a^1&\dots&a^k\end{bmatrix}^{(k)}=0$.

In the particular case~$k=n$, Eq.~\eqref{eq:volupp}
gives the well-known expression
\[
\text{volume} (\mathcal P(a^1,\dots,a^n))=|\det\begin{bmatrix} a^1&\dots& a^n \end{bmatrix} |.
\]

\begin{figure}
    \centering
    \begin{tikzpicture}
        \coordinate (a1) at (3,-0.5,0);
        \coordinate (a2) at (2,2,0);
        \coordinate (a3) at (1,0,-1);
        
        \draw[dashed,fill opacity=0.5,fill=green!10] (0,0,0) -- (a1) -- ($(a1) + (a2)$) -- (a2) -- cycle;
        \draw[dashed,fill opacity=0.5,fill=green!10] (a3) -- ($(a3) + (a1)$) -- ($(a3) + (a1) + (a2)$) -- ($(a3) + (a2)$) -- cycle;

        \draw[dashed,fill opacity=0.5,fill=green!10] (0,0,0) -- (a2) -- ($(a2) + (a3)$) -- (a3) -- cycle;
        \draw[dashed,fill opacity=0.5,fill=green!10] (a1) -- ($(a1) + (a2)$) -- ($(a1) + (a2) + (a3)$) -- ($(a1) + (a3)$) -- cycle;
        
        \draw[thick,->] (0,0,0)--(a1) node[right]{$a^1$};
        \draw[thick,->] (0,0,0)--(a2) node[above]{$a^2$};
        \draw[thick,->] (0,0,0)--(a3) node[above left=-0.1cm]{$a^3$};
        \draw (0,0,0) node[below]{0};
        \draw (3.15,1,0.5) node[]{$\mathcal{P}(a^1,a^2,a^3)$};

        \draw[->] (7,0,0) -- +(1,0,0) node[right]{$x$};
        \draw[->] (7,0,0) -- +(0,1,0) node[right]{$y$};
        \draw[->] (7,0,0) -- +(0,0,1) node[right]{$z$}; 
    \end{tikzpicture}
    \caption{A 3D parallelotope in~$\R^3$.}
    \label{fig:paralellotope}
\end{figure}

To study the evolution of a parallelotope
under a dynamical system, consider the linear time-varying~(LTV) system
\be\label{eq:ltv}
\dot x(t)=A(t)x(t),
\ee
with~$A:\R\to\R^{n\times n}$, 
and let~$x(t,x_0)$ denote its solution at time~$t$ when~$x(t_0)=x_0$. From hereon  we assume that the initial time~$t_0$ is zero. Then~$x(t,x_0)=\Phi(t )x_0$, where~$\Phi(t )\in\R^{n\times n} $ is the transition matrix from time~$ 0$ to time~$t$ of~\eqref{eq:ltv}. 
In other words, $\dot \Phi(t)=A(t)\Phi(t)$, with~$\Phi(0)=I_n$.

Fix initial conditions~$a^1,\dots,a^k\in\R^n$, and consider
the time-varying 
parallelotope~$
\mathcal P(x(t,a^1),\dots,x(t,a^k)).
$
%%%%
We already know that the volume of this 
parallelotope at time~$t$ is the~$L_2$ norm of the  vector 
$
v(t):=\begin{bmatrix} x(t,a^1)&\dots& x(t,a^k) \end{bmatrix} ^{(k)}
$. Now,
\begin{align}\label{eq:ddtvt}
%%%%%%%%%%%%%%%%%%
    \dot v(t)&= \frac{d}{dt} \big (\Phi(t)  \begin{bmatrix} a^1&\dots& a^k \end{bmatrix} \big )^{(k)}\nonumber\\
    &=\frac{d}{dt} \left ( (\Phi(t)) ^{(k)}\begin{bmatrix} a^1&\dots& a^k \end{bmatrix} ^{(k)}\right )
    \nonumber\\&=\frac{d}{dt} \left (( \Phi(t) )^{(k)}\right ) v(0)  ,
%%%%%%%%
\end{align}
%%%
where the second equality follows from the Cauchy-Binet theorem.

Eq.~\eqref{eq:ddtvt} naturally leads  to the notion of the $k$-additive compound of a (square) matrix~$B\in\R^{n\times n}$~\cite{fiedler_book}. For~$\varepsilon \geq  0$, consider the~$\binom{n}{k}\times\binom{n}{k}$ matrix~$(I_n+\varepsilon B)^{(k)}$. Every entry of this matrix is a~$k$-minor of~$I_n+\varepsilon B$, so 
%%%
\be\label{eq:expans}
(I_n+\varepsilon B)^{(k)}=L_0\varepsilon^0+L_1\varepsilon^1+\dots+L_k \varepsilon^k,
\ee
where every~$L_i$
is an~$r\times r$ matrix with~$r:=\binom{n}{k}$. 
%%%%
%%%%%%
Setting~$\varepsilon=0$ gives~$L_0=I_r$.
\begin{Definition}
    The matrix~$L_1$ in~\eqref{eq:expans} is called the~$k$-additive compound of~$B$, and is denoted~$B^{[k]}$.     
\end{Definition}

This definition implies that~$B^{[1]}=B$ and~$B^{[n]}=\tr(B)$. 
%%%%%%%%

Fix~$\varepsilon>0$. Then 
\begin{align*}
   \frac{ (\Phi(t+\varepsilon))^{(k)}- (\Phi(t ))^{(k)} }  {\varepsilon} &=
%%%
 \frac { ( \Phi(t)+\varepsilon A(t)\Phi(t) )^{(k)} +o(\varepsilon) -(\Phi(t ))^{(k)}    }{\varepsilon}\\
 &=
 \frac { (  I_n+\varepsilon A(t))^{(k)} (\Phi(t ))^{(k)}    +o(\varepsilon) -(\Phi(t ))^{(k)}    }{\varepsilon}\\
%%%
&=\frac{ (  I_n+\varepsilon A(t))^{(k)}+o(\varepsilon)-I_r}{\varepsilon}(\Phi(t ))^{(k)},
%%%%%
%%%
\end{align*}
and taking the limit~$\varepsilon\to 0 $  gives
\be\label{eq:expkk}
\frac{d}{dt} \left((\Phi(t) )^{(k)} \right) = (A(t))^{[k]} (\Phi(t) )^{(k)}.
\ee
Combining this with~\eqref{eq:ddtvt}, we conclude that
 \[
%%%%%%%%%%%%%%%%%%
    \frac{d}{dt} \left(\begin{bmatrix} x(t,a^1)&\dots& x(t,a^k) \end{bmatrix}   ^{(k)}\right) = (A(t))^{[k]} \begin{bmatrix} x(t,a^1)&\dots& x(t,a^k) \end{bmatrix} ^{(k)} .
%%%%%%%%
\]
%%%%%%%%%%%%%%
Roughly speaking, this implies that 
the evolution of~$k$-parallelotopes under the LTV dynamics~\eqref{eq:ltv} is   described by an LTV with the matrix~$(A(t))^{[k]}$. 
Note that when~\eqref{eq:ltv}
reduces to the LTI~$\dot x(t)=A x(t)$, Eq.~\eqref{eq:expkk} yields
\[
\frac{d}{dt} \left ( (\exp(At) )^{(k)} \right )=A^{[k]} (\exp(At) )^{(k)}.
\]

The term additive compound is justified by the fact that 
  for any~$A,B\in\R^{n\times n}$, we have~$(A + B)^{[k]} = A^{[k]} + B^{[k]}$.
%%%%%%%%%
Furthermore, it follows from the definitions above  that~$(A^{[k]})^T = (A^T)^{[k]}$, and that if~$T\in\R^{n\times n}$ is non-singular then
\be\label{eq:TATM1}
(TAT^{-1})^{[k]} = T^{(k)} A^{[k]}  (T^{(k)})^{-1}.
\ee
%%%%%

Using the definition of the $k$-additive compound and the Cauchy-Binet theorem yields  
\be\label{eq:effect_k_add}
A^{[k]} \begin{bmatrix} x^1&\dots&x^k \end{bmatrix}^{(k)} =
\begin{bmatrix} Ax^1&x^2&\dots&x^{k-1}&x^k \end{bmatrix}^{(k)}+\dots+
\begin{bmatrix} x^1&x^2&\dots &x^{k-1}&Ax^k \end{bmatrix}^{(k)}.
\ee

The~$k$-additive compound  satisfies an ``additive''  spectral property\updt{.}
If~$A\in\R^{n\times n}$  has eigenvalues~$\lambda_1,\dots,\lambda_n\in\C$ 
then the eigenvalues of~$A^{[k]}$ are all the~$\binom{n}{k}$ sums:
\[
\lambda_{i_1}+\dots+\lambda_{i_k}, \quad 1\leq i_1<\dots< i_k\leq n. 
\]
In particular,~$A^{[k]}$ is a  Hurwitz matrix  iff the sum of every $k$ eigenvalues of~$A $ has a negative real part, and then under the dynamics~$\dot x(t)=Ax(t)$ 
the volume of every~$k$-parallelotope decays to zero exponentially.

\subsection{Kronecker product }
%%%%%%%%%%%%
The Kronecker product of two matrices~$A\in\R^{n\times m}$ and~$B\in\R^{p\times q}$ is
\begin{equation}
    A \otimes B := \begin{bmatrix}
        a_{11} B &  \dots & a_{1m} B \\
        \vdots & \ddots & \vdots \\
        a_{n1} B & \dots & a_{nm} B
    \end{bmatrix}.
\end{equation}
Note that~$A \otimes B \in\R^{(np)\times(mq)}$.
We list several  properties of the   Kronecker product that will be used below  (see~\cite[Ch. 4]{Horn1991TopicsMatrixAna}
for proofs and more details). 
The product is distributive, that is,
 \begin{align*}
     A\otimes(B+C)&=A\otimes B+A\otimes C,\\
     (B+C)\otimes A&=B\otimes A+C\otimes A,
%%%
 \end{align*} 
 and associative, that is,
 \[
 (A\otimes B)\otimes C= A\otimes ( B \otimes C).
 \]
 Transposition is  distributive over the Kronecker product, that is,
 \[
 (A\otimes B)^T=A^T\otimes B^T. 
 \]
 %%%
The Kronecker product also satisfies  the following  mixed-product property: if~$A,B,C,D$ are matrices  such that the products~$AC$ and~$BD$ are well-defined then
\begin{equation}\label{eq:mixed_prod}
    (A \otimes B)(C \otimes D) = (AC) \otimes (BD).
\end{equation}

The Kronecker product is useful    when studying linear equations with matrix variables and, in particular, in  rewriting them as equations with vector variables. In particular, for $A  \in \R^{n \times n}$,
$B  \in \R^{m \times m}$,
and~$C,X \in\R^{n\times m}$,
  the matrix equation 
\begin{equation*}
    AX + XB  = C
\end{equation*}
is equivalent to the equation 
\begin{equation}\label{eq:vecxvecc}
  \left  ( (I_m \otimes  A)  +(B^T \otimes  I_n)\right )\vect(X) = \vect(C).
%%%%%%%%%%%%%%%%%%%%
\end{equation}

The Kronecker product also has a  useful  property with respect to the matrix exponential. If~$A_1,A_2 \in \R^{n\times n}$ then~\eqref{eq:mixed_prod} implies that the matrices~$A_1\otimes I_n$ and~$I_n\otimes A_2$ commute, so 
\begin{align}\label{eq:exp2}
    \exp(A_1\otimes I_n + I_n\otimes A_2) &=\exp(A_1\otimes I_n)\exp( I_n\otimes A_2)\nonumber\\
    &=(\exp(A_1)\otimes I_n  )( I_n\otimes \exp(A_2) ) \nonumber\\
    &=\exp(A_1)\otimes \exp(A_2),
\end{align}
where the last equation follows from using~\eqref{eq:mixed_prod} again.

More generally, for
$A_1,\dots,A_s \in \R^{n \times n}$, we have
\begin{equation}\label{eq:kron_exp}
%%%%%%%%%%%%%%%%
\exp\left( A_1 \otimes I_n \otimes \dots \otimes I_n + I_n \otimes A_2 \otimes \dots \otimes I_n + I_n \otimes \dots \otimes I_n\otimes  A_s\right) =  \exp(A_1) \otimes \dots\otimes \exp(A_s),
%%%%%%%%%%%%%%%%%%%%%%%%
\end{equation}
where every Kronecker product includes $s$ terms. 

 The  ``$k$th Kronecker power'', denoted~$A^{\otimes k}$,  
is defined inductively
by
\[
A^{\otimes k}:=  A^{\otimes (k-1)} \otimes A,
\]
with~$A^{\otimes 1}   :=A$. 
If~$A\in\R^{n\times n}$ with eigenvalues~$\lambda_1,\dots,\lambda_n$ then the eigenvalues of~$A^{\otimes k}$ are all the~$n^k$ products
\[
\lambda_{i_1}\dots \lambda_{i_k}, \text{ with } i_\ell \in \{1,\dots,n \}.
\]
Note that this set includes in particular all the eigenvalues of~$A^{(k)}$.

The next example demonstrates 
how the Kronecker product naturally  appears   
in the multiplicative 
compound  of  a block matrix. 

%%%%%%%%%%%%%%%%%%%%%%%%%%%%%%%%%%%%%%%%%%
 \begin{Example}\label{exa:multi_mixed}
 %%%%%%%%%%%%%%%%%%%%%%%%%%%%%%%%%%%%%%%%%
    Consider the matrix~$B=\diag(d_1,\dots,d_4)$. 
    Then a direct calculation gives
    \[
      B ^{(2)}=\diag(d_1d_2,d_1d_3,d_1 d_4, d_2 d_3,d_2 d_4,d_3 d_4 ).
    \]
Note that if we write~$B$ as a block matrix~$B=\diag(B_1,B_2)$,
with~$B_1=\diag(d_1,d_2)$ and~$B_2=\diag(d_3,d_4)$,  then 
\[
B^{(2)} = 
   \diag( B_1^{(2)} , B_1\otimes B_2,B_2^{(2)}).
\]
%%%%
%%%%%
\end{Example}

\subsection{The Kronecker sum}
%%%%
The Kronecker sum of two square matrices~$A\in\R^{n\times n},B\in\R^{m\times m}$ is the~$(nm)\times(nm)$ matrix
\[
A\oplus B := A\otimes I_m+I_n\otimes B. 
\]
With this notation,~\eqref{eq:exp2} becomes  
\[
\exp(A_1\oplus A_2)=\exp(A_1)\otimes\exp(A_2),
\]
and~\eqref{eq:vecxvecc} becomes
\be\label{eq:btvec}
(B^T\oplus A)\vect(X)=\vect(C). 
\ee

The next example demonstrates 
how the Kronecker sum naturally appears   
in the additive compound  of  a block matrix. 

%%%%%%%%%%%%%%%%%%%%%%%%%%%%%%%%%%%%%%%%%%
 \begin{Example}\label{exa:mixed}
 %%%%%%%%%%%%%%%%%%%%%%%%%%%%%%%%%%%%%%%%%
    Consider the matrix~$B=\diag(d_1,\dots,d_4)$. 
    Then a direct calculation gives
    \[
    (I_4+\varepsilon B )^{(2)}=\diag(s_1 s_2,s_1 s_3,s_1 s_4,s_2 s_3, s_2 s_4,s_3 s_4), 
    \]
    where~$s_i:=1+\varepsilon d_i$, and thus~\eqref{eq:expans} with~$n=4$ and~$k=2$ gives 
    \[  
B^{[2]}= \diag(d_1+d_2,d_1+d_3,d_1+d_4,d_2+d_3,d_2+d_4,d_3+d_4). 
\]
Writing~$B$ as a block matrix~$B=\diag(B_1,B_2)$,
with~$B_1=\diag(d_1,d_2)$ and~$B_2=\diag(d_3,d_4)$ yields  
\[
B^{[2]} = 
   \diag( B_1^{[2]} , B_1\oplus B_2,B_2^{[2]}).
\]
%%%%
%%%%%
\end{Example}
 
For~$A\in\R^{n\times n}$ 
the  ``$k$th Kronecker sum'' of~$A$ is defined   by
\be\label{eq:def_k_kron_sum}
A^{\oplus k}:=   A\otimes I_n \otimes \dots \otimes I_n+ 
   I_n\otimes A  \otimes I_n   \otimes\dots \otimes I_n+
\dots +
  I_n \otimes \dots \otimes I_n \otimes A, 
\ee
where every product includes~$k$ terms. 
For example,~$A^{\oplus 1}=A$, and~$A^{\oplus 2}= A\oplus A$.

If~$A\in\R^{n\times n}$ with eigenvalues~$\lambda_1,\dots,\lambda_n$ then the eigenvalues of~$A^{\oplus k}$ are all the~$n^k$ sums
\[
\lambda_{i_1}+\dots+ \lambda_{i_k}, \text{ with } i_\ell \in \{1,\dots,k \}.
\]
Note that this set includes in particular all the eigenvalues of~$A^{[k]}$.

%%%%%%%%%%%%%
\section{Main results}\label{sec:main}
%%%%%%%%%%%%%%%%%%%%%%%%%%%%%%
Fix   integers~$n>1$ and~$k\in\{1,\dots,n\}$. 
Let~$r:=\binom{n}{k}$. 
 Fix~$x^1,\dots,x^k \in \R^n$. 
 Let~$X:=\begin{bmatrix}
     x^1&\dots&x^k 
 \end{bmatrix}\in\R^{n\times k}$, 
 and let~$z:=\begin{bmatrix}
     x^1&\dots&x^k 
 \end{bmatrix}^{(k)}\in\R^{r\times 1}$.
 The vector~$z$ includes all
 the~$k\times k$ minors of~$\begin{bmatrix}
     x^1&\dots& x^k 
 \end{bmatrix}$.
Our first goal is to relate this to  Kronecker products of the~$x^i$s. This    requires defining  two auxiliary matrices that we denote by~$M_{n,k}$ and~$L_{n,k}$.  

\subsection{The matrices~$M_{n,k}$ and $L_{n,k}$}
%%%%%%%%%%%%%%%%%%%%%%%%%%%%%%%%%%

We   begin by  defining 
the matrices~$M_{n,k} $ and~$L_{n,k} $.
%%%%%%%%%%%%%%%%%%%%%
\begin{Definition}\label{def:mnln}
Let~$e^1,\dots,e^n$ denote the canonical  basis of~$\R^n$. Fix~$k\in\{1,\dots,n\}$, and define matrices  
%%%%
\[
M_{n,k}\in\R^{n^k \times \binom{n}{k} }  \text{ and } L_{n,k}\in\R^{ \binom{n}{k}  \times n^k}
\]
by
%%%
\begin{align}\label{eq:mn_alternative}
    	M_{n,k} &:= \sum_{(i_1,\dots,i_k) \in Q(n,k)} \sum_{(j_1,\dots,j_k) \in S(i_1,\dots,i_k)} \sigma(j_1,\dots,j_k) (e^{j_1} \otimes \dots \otimes e^{j_k}) \left(\begin{bmatrix} e^{i_1} & \dots & e^{i_k} \end{bmatrix}^{(k)}\right)^T, \nonumber\\
    %%%%%%%
    	L_{n,k} &:= \sum_{(i_1,\dots,i_k) \in Q(n,k)} \begin{bmatrix} e^{i_1} & \dots & e^{i_k} \end{bmatrix}^{(k)} (e^{i_1} \otimes \dots \otimes e^{i_k})^T,
    \end{align}
    where $S(i_1,\dots,i_k)$ is  the set of all permutations of the sequence $i_1,\dots,i_k$, and $\sigma(j_1,\dots,j_k)\in\{-1,1\}$ is  the signature of   the permutation~$j_1,\dots,j_k$.
\end{Definition}

\begin{Example}\label{exa:mnln}
    %%%%
For~$k=n$ the only element in~$Q(n,k)$ is the sequence~$(1,2,\dots,n)$, so
    %%%
    \begin{align*}
    %%%
    	M_{n,n} &=   \sum_{(j_1,\dots,j_n) \in S(1,\dots,n)} \sigma(j_1,\dots,j_n) (e^{j_1} \otimes \dots \otimes e^{j_n}) (\begin{bmatrix} e^{ 1} & \dots & e^{ n} \end{bmatrix}^{(n)})^T   \\
     &=  \sum_{(j_1,\dots,j_n) \in S(1,\dots,n)} \sigma(j_1,\dots,j_n) (e^{j_1} \otimes \dots \otimes e^{j_n}) , 
    %%%%%%%
     \end{align*}
    and
    \begin{align*}
    	L_{n,n} &=  \begin{bmatrix} e^{ 1} & \dots & e^{ n} \end{bmatrix}^{(n)}  (e^{1} \otimes \dots \otimes e^{n})^T\\
     &=(e^1)^T\otimes\dots\otimes(e^n)^T .
    \end{align*}
    \end{Example}

\subsection{Properties of the matrices~$M_{n,k}$ and~$L_{n,k}$}
%%%%%%%%%%%%%%%%%%%%%%%%%
The matrices~$M_{n,k}$ and~$L_{n,k}$ play an important role in this paper, so we begin by analyzing some of their properties. 
It is useful to rewrite~$M_{n,k}$  and~$L_{n,k}$ in another form. 
To do this 
   requires the following definition. 

%%%%%%%%%%%%     
\begin{Definition}
  %%%%%%%%%%%%%%%
    Let~$w^1(n,k) ,\dots,w^r(n,k)$ denote the canonical   basis of~$\R^r$. Let~$v^1(n,k) ,\dots,v^{n^k}(n,k)$ denote the canonical  basis of~$\R^{n^k}$.
\end{Definition}

 %%%%%%%   
To simplify the notation, we will usually just write~$w^i$ and~$v^i$, omitting the dependence on~$n,k$, but it is important to remember  that the dimension of these vectors depends on~$n,k$. 
%%%%%%%%%%%%%%%%%%%
 For example, for~$n=4$ and~$k=2$, we have that~$r=\binom{4}{2}=6$,  $w^1,\dots,w^{6}$
denote the canonical basis of~$\R^{6}$, $n^k=16$, and~$v^1,\dots,v^{16}$
 denote the canonical basis of~$\R^{16}$.

\begin{Lemma}\label{lem:alternative}
%%%%%%%%%%%%%%%%%%%%%%%%%
The matrices~$M_{n,k}$ and~$L_{n,k}$ satisfy 
\begin{align}
	M_{n,k} & = \sum_{(i_1,\dots,i_k) \in Q(n,k)} \sum_{(j_1,\dots,j_k) \in S(i_1,\dots,i_k)} \sigma(j_1,\dots,j_k) v^{\ind^R_{n,k}(j_1,\dots,j_k)} (w^{\ind^Q_{n,k}(i_1,\dots,i_k)})^T, \\
%%%%%%%
	L_{n,k} &  = \sum_{(i_1,\dots,i_k) \in Q(n,k)} w^{\ind^Q_{n,k}(i_1,\dots,i_k)}(v^{\ind^R_{n,k}(i_1,\dots,i_k)})^T . 
\end{align}
%%%%%%%%%
\end{Lemma}

\begin{IEEEproof}
    %%%%
 It is shown in
  Lemma~\ref{lemma:index_R} in the Appendix that 
  \[    \ind^R_{n,k}(i_1,\dots,i_k) = \sum_{\ell=1}^k (i_\ell - 1)n^{k-\ell} + 1.
  \]
Combining this with the definition of the Kronecker product implies that~$(x^1 \otimes \dots \otimes x^k)_{\ind^R_{n,k}(i_1,\dots,i_k)} = x^1_{i_1} \dotsm x^k_{i_k}$. 
In particular,   for the canonical basis   
  $e^1,\dots,e^n$   of~$\R^n$, we have 
  $(e^{i_1} \otimes \dots \otimes e^{i_k})_{\ind^R_{n,k}(i_1,\dots,i_k)} = e^{i_1}_{i_1} \dotsm e^{i_k}_{i_k}=1$,  so
%%%%
\begin{equation}      
e^{i_1} \otimes \dots \otimes e^{i_k}=v^{\ind^R_{n,k}(i_1,\dots,i_k)} .
    \end{equation}
    
    Similarly, by the definition of the multiplicative compound, for any~$(i_1,\dots,i_k) \in Q(n,k)$, we have that $\begin{bmatrix} e^{i_1} & \dots & e^{i_k} \end{bmatrix}^{(k)} = w^{\ind^Q_{n,k}(i_1,\dots,i_k)}$.
    Substituting these expressions in Definition~\ref{def:mnln} completes the proof. 
%%%%%%%%%%%%
%%%%%%%%%%%%%
\end{IEEEproof}

\begin{Lemma}\label{lemma:LM_IS_I}
    The matrices~$M_{n,k}$ and~$L_{n,k}$ satisfy
       \be\label{eq:LM_is_Identity}
    L_{n,k}M_{n,k}=I_r. 
    \ee 
\end{Lemma}
%%%%
\begin{IEEEproof}
    %%%%%%%%%%%%%%%%
 We begin by analyzing~$\image(M_{n,k})$. 
 Fix an index~$p\in \{1,\dots,r\}$,
    and consider
    \begin{align*}
       y:&= M_{n,k}   w^p\\ &= \sum_{(i_1,\dots,i_k) \in Q(n,k)} \sum_{(j_1,\dots,j_k) \in S(i_1,\dots,i_k)} \sigma(j_1,\dots,j_k) v^{\ind^R_{n,k}(j_1,\dots,j_k)} (w^{\ind^Q_{n,k}(i_1,\dots,i_k)})^T w^p.
    \end{align*}
    Let~$(i_1',\dots,i_k')$ be the sequence in~$Q(n,k)$ such that~$ \ind^Q_{n,k}(i_1',\dots,i_k') =p$. Then 
     \begin{align*}
       y&=  \sum_{(j_1,\dots,j_k) \in S(i_1',\dots,i_k')} \sigma(j_1,\dots,j_k) v^{\ind^R_{n,k}(j_1,\dots,j_k)}  .
    \end{align*}
%%%%
Thus, the multiplication by~$M_{n,k}$ ``duplicates''
the single one  in  entry~$p$ of~$w^p$ 
to all the~$k!$  entries with indices~$\ind^R_{n,k} (j_1,\dots,j_k) $,  $(j_1,\dots,j_k) \in S(i_1',\dots,i_k')$,   in~$M_{n,k}w^p$, each 
with either a plus or a minus sign (see Fig.~\ref{fig:LnkMnk}). 
Now, 
   \begin{align*}
    L_{n,k}   y&=  \sum_{(i_1,\dots,i_k) \in Q(n,k)} w^{\ind^Q_{n,k}(i_1,\dots,i_k)}(v^{\ind^R_{n,k}(i_1,\dots,i_k)})^T \sum_{(j_1,\dots,j_k) \in S(i_1',\dots,i_k')} \sigma(j_1,\dots,j_k) v^{\ind^R_{n,k}(j_1,\dots,j_k)}  \\
    &= w^{\ind^Q_{n,k}(i_1',\dots,i_k')}   \sigma(i'_1,\dots,i'_k) \\
    &=w^p. 
    %%%%%
    \end{align*}
    Since the index~$p$ is arbitrary, this proves~\eqref{eq:LM_is_Identity}. 
%%%%%%%%%
\end{IEEEproof}

\begin{figure} 
    \centering
    \begin{equation*}
        \begin{bmatrix}
            0 \\
            a \\
            b \\
            -a \\
            0 \\
            c \\
            -b \\
            -c \\
            0 \\
        \end{bmatrix} = M_{3,2} \begin{bmatrix}
            a \\
            b \\
            c
        \end{bmatrix}, \text{ and }
        \begin{bmatrix}
            b \\
            c \\
            f
        \end{bmatrix} = L_{3,2} \begin{bmatrix}
            a \\
            b \\
            c \\
            d \\
            e \\
            f \\
            g \\
            h \\
            i \\
        \end{bmatrix}.
    \end{equation*} \caption{Multiplication by the matrices matrices~$M_{3,2}$ and $L_{3,2}$. 
        Note in particular that $
             L_{3,2} M_{3,2} \begin{bmatrix}
                a \\
                b \\
                c
            \end{bmatrix}= \begin{bmatrix}
                a \\
                b \\
                c
            \end{bmatrix}.
        $\label{fig:LnkMnk}}
%%%%
\end{figure}

\begin{Example}
 %%%%%%%%%%  
Consider the case~$k=1$. Then
\begin{align*}
	M_{n,1} & = \sum_{i \in Q(n,1)}   v^{\ind^R_{n,1}(i)} (w^{\ind^Q_{n,1}(i )})^T, \\
 &=\sum_{i =1}^n   v^{i} (w^i)^T\\
 &=I_n. 
%%%%%%%
\end{align*}
 
 Now   consider the case~$k=2$. Then
    \begin{align}\label{eq:M_n,2}
%%%%%%%%%%%%%%%%%%%    
	M_{n,2} & = \sum_{(i_1,i_2) \in Q(n,2)} \sum_{(j_1,j_2) \in S(i_1,i_2)} \sigma(j_1,j_2) v^{\ind^R_{n,2}(j_1,j_2)} (w^{\ind^Q_{n,2}(i_1,i_2)})^T\nonumber \\
 %%%
 &= \sum_{(i_1,i_2) \in Q(n,2)} \left(
 %%%%%%%
  v^{(i_1-1)n+i_2} (w^{\ind^Q_{n,2}(i_1,i_2)})^T
-v^{(i_2-1)n+i_1} (w^{\ind^Q_{n,2}(i_1,i_2)})^T
 \right ).
 %%%%%%%%%%%%%%%%
\end{align}
%%%
In particular, if~$n=2$ then~$Q(n,2)=Q(2,2)$ includes only the sequence~$(i_1,i_2)=(1,2)$, so
%%%%%%%%%%%%%%%%%%%%%%%%%
    \begin{align*} 
%%% 
	M_{2,2} 
 %%%
 &=  
 %%%%%%%
  v^{2} (w^{1})^T
-v^{3} (w^{1})^T\\
    &=\begin{bmatrix}
        0& 1&-1 & 0 
    \end{bmatrix}^T. 
\end{align*}
%%%%%%%%%
 %%%%%%   
 If~$n=3$ then~$Q(n,2)=Q(3,2)=((1,2),(1,3),(2,3))$, so
    \begin{align}\label{eq:M_3,2}
%%%%%%%%%%%%%%%%%%%    
	M_{3,2} 
 &= \sum_{(i_1,i_2) \in Q(n,2)} \left(
 %%%%%%%
  v^{3(i_1-1)+i_2} (w^{\ind^Q_{3,2}(i_1,i_2)})^T
-v^{3(i_2-1)+i_1} (w^{\ind^Q_{3,2}(i_1,i_2)})^T
 \right )\nonumber\\
 &=
 (v^2 -v^4)(w^1)^T+
 (v^3 -v^7)(w^2)^T+
 (v^6 -v^8)(w^3)^T\nonumber\\
 &=\begin{bmatrix}
     0&0&0\\
     1&0&0\\
     0&1&0\\
     -1&0&0\\
     0&0&0\\
     0&0&1\\
     0&-1&0\\
     0&0&-1\\
     0&0&0
 \end{bmatrix}
 .
 %%%%%%%%%%%%%%%%
\end{align}

%%%%%%%%%%%%%%%%

 Similarly, $L_{n,2}  \in \R^{\binom{n}{2} \times n^2}$ is given by 
\begin{align}
%%%%%%%%%%%%%%%%%%
        L_{n,2} &= 	\sum_{(i_1,i_2) \in Q(n,2)} w^{\ind^Q_{n,2}(i_1,i_2)}(v^ {  (i_1-1)n +(i_2-1)+1   })^T,
%%%%%%%%%%%%%%%%%
\end{align}
In particular, if~$n=2$ then
\begin{align*}
    L_{2,2} &= w^1 (v^2)^T  \\
    &=\begin{bmatrix}
        0& 1&0 & 0 
    \end{bmatrix}. 
\end{align*}
%%%%%%%%%%%%%%%%%%
Note that
\begin{align*}
    L_{2,2}M_{2,2}&=\begin{bmatrix}
        0& 1&0 & 0 
    \end{bmatrix}\begin{bmatrix}
        0& 1&-1 & 0 
    \end{bmatrix}^T \\
    &=I_1,
\end{align*}
(see~\eqref{eq:LM_is_Identity}). 
%%%%%%%%%
Now consider the case~$k=n$. Then~$r=\binom{n}{n}=1$, so~$w^1=1$ and~$Q(n,n) $ includes only the sequence~$(1,\dots,n)$,
so 
\begin{align}\label{eq:mnn}
	M_{n,n} &=   \sum_{(j_1,\dots,j_k) \in S(1,\dots,n)} \sigma(j_1,\dots,j_k) v^{\ind^R_{n,n}(j_1,\dots,j_k)} ,  
\end{align}
and
\begin{equation}\label{eq:lnn}
L_{n,n}  =   (v^{\ind^R_{n,n}(1,\dots,n)})^T 
     =  (v^{\frac{n^n - n}{(n-1)^2}})^T , 
\end{equation}
%%%
\end{Example}
where the last equality is proven in the Appendix.

\begin{Remark}\label{rem:Ln_Mn_vec}
%%%%%%%%%%%%%
%%%%%%%%%%%%%%
Recall that a matrix~$A\in\R ^{ n \times n}$ is called skew-symmetric if~$A+A^T=0$. The half-vectorization of such a matrix, denoted~$\vech(A)$,  is the column vector obtained by stacking the entries below    the main diagonal of~$A$ column by column.
For example, for the~$3\times 3$ skew-symmetric matrix
\be\label{eq:A33_skew}
A_3=\begin{bmatrix}
0& a_{12}&a_{13}\\
-a_{12}& 0&a_{23}  \\
-a_{13} & -a_{23} &0  
\end{bmatrix},
\ee 
we have 
\[
\vech( A_3 ) =\begin{bmatrix}
-    a_{12} &-a_{13}& -a_{23} 
\end{bmatrix}^T. 
\]
%%%%%%%
For an~$n\times n$ skew-symmetric matrix~$A$,
    the matrices $M_{n,2}$ and $L_{n,2}$    relate between vectorizations and half-vectorizations of~$A$:
\be\label{eq:vec_and_vech}
    \vect(A) = M_{n,2} \vech(A) , \text{ and } \vech(A) = L_{n,2} \vect(A).
\ee
  For example, for the case~$n=3$  using~$M_{3,2}$ in~\eqref{eq:M_3,2} gives 
%%%%%%%%%%%%%%%%%%%%%%%
   \begin{align*}
        M_{3,2}\vech(A_3)&= M_{3,2}\begin{bmatrix}
    -a_{12} &-a_{13}& -a_{23} 
\end{bmatrix}^T\\
&=\begin{bmatrix}
    0&-a_{12}&-a_{13}&a_{12}&0&-a_{23}&a_{13}&a_{23} &0
\end{bmatrix}^T\\&=\vect(A_3). 
%%%%      
\end{align*}
%%%%%%  
    The  relations~\eqref{eq:vec_and_vech}
    have been used in~\cite{Angeli2023smallgain} to study skew-symmetric solutions of matrix ODEs of the form~$\dot X = AX + XA^T$, where $A$ and $X$ are block matrices.
\end{Remark}
%%%%%%

\begin{Remark}
    Similar matrices, sometimes referred to  as duplication and elimination matrices, also appear when considering vectorizations and half-vectorizations of \emph{symmetric} matrices, see for example~\cite{Magnus1980Elimination}.
\end{Remark}

 \subsection{Algebraic expressions for the $k$-multiplicative compound}
%%%%%%%
We can now state one of the  main results in this section. This result allows to express the $k$-multiplicative compound of an~$n\times k$ matrix using  a sum of Kronecker products.  The matrices~$L_{n,k}$ and~$M_{n,k}$ link the   two representations. 
%%%%%%
\begin{Theorem}
\label{thm:MNK+and_LNK}
%%%%%%%%%%%%%%%%%%%
    For any~$x^1,\dots,x^k \in \R^n$ with $k \le n$, we have 
%%%
\begin{align}
%%%%%%%%%%%%%%%%%%%
\begin{bmatrix}x^1 & \dots & x^k\end{bmatrix}^{(k)}& = L_{n,k} \sum_{(j_1,\dots,j_k) \in S(1,\dots,k)} \sigma(j_1,\dots,j_k) (x^{j_1} \otimes \dots \otimes x^{j_k}) , \label{eq:lnk_general}
\end{align}
and
\begin{align}
%%%%%
    M_{n,k} \begin{bmatrix}x^1 & \dots & x^k\end{bmatrix}^{(k)}& = \sum_{(j_1,\dots,j_k) \in S(1,\dots,k)} \sigma(j_1,\dots,j_k) (x^{j_1} \otimes \dots \otimes x^{j_k}) .\label{eq:mnk_general}
%%%%%%%%%%
\end{align}
%%%%%%%%%%%%%%%%%%%%%%%%%%%%%
\end{Theorem}

\begin{Example}
    For~$k =2$ Theorem~\ref{thm:MNK+and_LNK} becomes   
%%%%%%%%%%%%%%%%%%%
\begin{align*}
    \begin{bmatrix}x^1 & x^2\end{bmatrix}^{(2)}& = L_{n,2} (x^1 \otimes x^2 - x^2 \otimes x^1) \label{eq:kron_to_mul_comp}   , 
    \end{align*}
and
\begin{align*}
%%%%%%
    M_{n,2} \begin{bmatrix}x^1 & x^2\end{bmatrix}^{(2)}& =  x^1 \otimes x^2 - x^2 \otimes x^1   .
\end{align*}
%%%%%%%%%%%%%%%%%%%%%%%%%%%%%
 Intuitively speaking, this can be explained as follows. 
The vector~$\begin{bmatrix}x^1 & x^2\end{bmatrix}^{(2)}\in\R^{ \binom{n}{2}}$ is the set of    all the two~$2\times 2$ minors of~$\begin{bmatrix}x^1 & x^2\end{bmatrix}$. The entries of the vector~$(x^1 \otimes x^2 - x^2 \otimes x^1)\in\R^{n^2} $ also include all these minors,  and multiplication by
the matrices~$L_{n,2}$ and~$M_{n,2}$       transforms one vector into the other.
%%%%%%%%%
 \end{Example}

\begin{IEEEproof}[Proof of Theorem~\ref{thm:MNK+and_LNK}]
%%%%%
 Let 
 \[
 z:= \sum_{(j_1,\dots,j_k) \in S(1,\dots,k)} \sigma(j_1,\dots,j_k) (x^{j_1} \otimes \dots \otimes x^{j_k}).
 \]
 We begin by considering the particular  case  
 the   where  every~$x^i$ is a canonical  vector in~$\R^n$, that is, $\begin{bmatrix}
     x^1&\dots&x^k
 \end{bmatrix}=
 \begin{bmatrix}e^{i_1}& \dots&e^{i_k} 
 \end{bmatrix}$. 
 We consider three  sub-cases.
 
 \emph{Case 1.}
 %%%%%%%%%%%%%
Suppose that  there exist two indices~$i_p$, $i_q$ with~$p\not =q$ and~$i_p=i_q$.  Then~$ \begin{bmatrix}x^1 & \dots & x^k\end{bmatrix}^{(k)}=0$.
 %%%%%
 Any term~$x^{j_1} \otimes \dots \otimes x^{j_k}$ includes~$x^{i_p},x^{i_q}$, and permuting between these two vectors gives another term in the sum~$z$, but with an opposite sign. Thus,~$z=0$.  So in this case, equations~\eqref{eq:lnk_general} and~\eqref{eq:mnk_general} hold.  
 
 \emph{Case 2.}
 %%%%%%%%%%%%%%%%
 Suppose that~$x^1=e^{i_1},\dots,x^k=e^{i_k}$, with
 \be\label{eq:i1ikorder}
 i_1<\dots<i_k.
 \ee
 Then
\[ 
    \begin{bmatrix}x^1 & \dots & x^k\end{bmatrix}^{(k)}= w^{\ind^Q_{n,k}(i_1,\dots,i_k)}.
 \]
Consider
\begin{align*}
L_{n,k}z & =  \sum_{(i_1,\dots,i_k) \in Q(n,k)} w^{\ind^Q_{n,k}(i_1,\dots,i_k)}(v^{\ind^R_{n,k}(i_1,\dots,i_k)})^T  \sum_{(j_1,\dots,j_k) \in S(1,\dots,k)} \sigma(j_1,\dots,j_k) v^{\ind^R_{n,k}(i_{j_1},\dots,i_{j_k})}\\
%%%%
&= \sigma(i_1,\dots,i_k) w^{\ind^Q_{n,k}(i_1,\dots,i_k)}\\
&=w^{\ind^Q_{n,k}(i_1,\dots,i_k)},
\end{align*}
where the last step follows from~\eqref{eq:i1ikorder}. 
Thus,~\eqref{eq:lnk_general} holds. To verify~\eqref{eq:mnk_general},   note that~$  M_{n,k} \begin{bmatrix}x^1 & \dots & x^k\end{bmatrix}^{(k)}$ is equal to
\begin{align*}
 %     M_{n,k} \begin{bmatrix}x^1 & \dots & x^k\end{bmatrix}^{(k)} & = 
& \sum_{(p_1,\dots,p_k) \in Q(n,k)} \sum_{(j_1,\dots,j_k) \in S(p_1,\dots,p_k)} \sigma(j_1,\dots,j_k) v^{\ind^R_{n,k}(j_1,\dots,j_k)} (w^{\ind^Q_{n,k}(p_1,\dots,p_k)})^T w^{\ind^Q_{n,k}(i_1,\dots,i_k)}\\
%%%%%
      &= \sum_{(j_1,\dots,j_k) \in S(i_1,\dots,i_k)} \sigma(j_1,\dots,j_k) v^{\ind^R_{n,k}(j_1,\dots,j_k)},
%%%% 
\end{align*}
implying  that~\eqref{eq:mnk_general} holds.  

 %%%%%%%%%%%%%%%%%%%%%%%%%%
\emph{Case 3.}
 %%%%%%%%%%%%%%%%%%%%%%%%%%
  Suppose that~$x^1=e^{i_1},\dots,x^k=e^{i_k}$, with~$i_j\not=i_\ell$ for all~$j\not = \ell$.
Then there exists a permutation~$\ell_1,\dots,\ell_k$ of~$i_1,\dots,i_k$ such that~$\ell_1<\dots<\ell_k$. 
Then~$\begin{bmatrix}x^1 & \dots & x^k\end{bmatrix}^{(k)}=\sigma(i_1,\dots,i_k)  \begin{bmatrix} e^{\ell_1} & \dots & e^{\ell_k}\end{bmatrix}^{(k)} $. Let~$s:=\sigma(i_1,\dots,i_k) $.
  By Case 2,
\begin{align*}
%%%%
       \begin{bmatrix}x^1 & \dots & x^k\end{bmatrix}^{(k)}&=s  \begin{bmatrix} e^{\ell_1} & \dots & e^{\ell_k}\end{bmatrix}^{(k)}\\
%%%
&= s L_{n,k} \sum_{(j_1,\dots,j_k) \in S(\ell_1,\dots,\ell_k)} \sigma(j_1,\dots,j_k) (x^{j_1} \otimes \dots \otimes x^{j_k}),
%%%%
\end{align*}
and this gives
\begin{align*}
%%%%
       \begin{bmatrix}x^1 & \dots & x^k\end{bmatrix}^{(k)}&= 
   L_{n,k} \sum_{(j_1,\dots,j_k) \in S(i_1,\dots,i_k)} \sigma(j_1,\dots,j_k)   (x^{j_1} \otimes \dots \otimes x^{j_k}),
%%%%
\end{align*}
so~\eqref{eq:lnk_general} holds. A similar argument shows that~\eqref{eq:mnk_general} also holds.
%%%
%%
  This  completes the proof of Theorem~\ref{thm:MNK+and_LNK}
     when every~$x^i$ is a  standard basis vector in~$\R^n$. For arbitrary vectors in~$\R^n$ the proof  follows from the multilinearity of the multiplicative compound and of the Kronecker product.
%%%%%%%%%%%%%%%%%%%%%%%%%%%%  
\end{IEEEproof}

 The next result describes a useful  relation between Kronecker products  and multiplicative compounds. 
\begin{Theorem}\label{thm:akk}
%%%%
For any~$A \in\R^{n\times m}$ and $k \in \{1,\dots,\min\{n,m\}\}$, we have  
\be\label{eq:aok}
A^{\otimes k} M_{m,k}= M_{n,k} A^{(k)}  . 
\ee
%%%%%
\end{Theorem}

\begin{IEEEproof}
%%%%
Pick  arbitrary~$x^1,\dots,x^k \in \R^m$,  and let~$y:=\begin{bmatrix} x^1&\dots&x^k \end{bmatrix}^{(k)}$. Then using~\eqref{eq:mnk_general} and the mixed product rule 
yields 
%%%%%%%%%%%
\begin{align*}
    A^{\otimes k} M_{m,k}  y & =  A^{\otimes k} \sum_{(j_1,\dots,j_k) \in S(1,\dots,k)} \sigma(j_1,\dots,j_k) (x^{j_1} \otimes \dots \otimes x^{j_k}) \\
    &= \sum_{(j_1,\dots,j_k) \in S(1,\dots,k)} \sigma(j_1,\dots,j_k) (Ax^{j_1} \otimes \dots \otimes Ax^{j_k}) . 
\end{align*}
Using~\eqref{eq:mnk_general}  again  and   the Cauchy-Binet theorem gives
\begin{align*}
    A^{\otimes k} M_{m,k}   y   
    &=  M_{n,k} \begin{bmatrix} Ax^1&\dots&Ax^k \end{bmatrix}^{(k)}  \\
    &=  M_{n,k} A^{(k)} \begin{bmatrix} x^1&\dots& x^k \end{bmatrix}^{(k)} \\
    &=M_{n,k}A^{(k)} y, 
\end{align*}
and since~$x^1,\dots,x^k$ are arbitrary,  this proves~\eqref{eq:aok}. 
%%%%%
\end{IEEEproof}

The next result uses the expressions derived  above  
to provide  a formula  for the~$k$-multiplicative compound of a    matrix~$A$ using   the $k$th Kronecker power of~$A$.

\begin{Theorem}\label{thm:multi_compi}
%%%%%%%%%%%%%%%%
 Let~$A\in\R^{n \times m}$,  and fix~$k \in \{1,\dots,\min\{n,m\}\}$. Then 
 \be
    \label{eq:mul_kron}
        A^{(k)} = L_{n,k}   A^{\otimes k} M_{m,k}  . 
    \ee
%%%
%%%%%%%%
\end{Theorem}

\begin{IEEEproof}
%%%%
Multiply both sides of~\eqref{eq:aok} by~$L_{n,k}$ from the left and apply~\eqref{eq:LM_is_Identity}.
%%%%%%%%%%%%%%%%%%%%%%%%
\end{IEEEproof}

For example, for~$k=2$,  Eq.~\eqref{eq:mul_kron} gives
 \[
        A^{(2)} = L_{n,2} \left( A \otimes A\right) M_{m,2}  . 
    \]

\begin{Remark}
    As a ``sanity check'', we show that the expression in Theorem~\ref{thm:multi_compi} implies   the Cauchy-Binet theorem. 
   Let~$A\in\R^{n\times m}$ and~$B\in\R^{m\times p}$, and fix~$k\in\{1,\dots,\min\{n,m,p\}\}$.   Theorem~\ref{thm:multi_compi} gives 
    \begin{align*}
        A^{(k)}B^{(k)}&=  L_{n,k}     A^{\otimes k}   M_{m,k}   L_{m,k}     B^{\otimes k}  M_{p,k}  .
    \end{align*}
    By Theorem~\ref{thm:akk}, $B^{\otimes k}  M_{p,k}=M_{m,k} B^{(k)}$, 
    and using this and~\eqref{eq:LM_is_Identity} gives  
 $   
        A^{(k)}B^{(k)} =  L_{n,k}     A^{\otimes k}   M_{m,k}    B^{(k)}  .
$     
   Using Theorem~\ref{thm:akk} again yields 
       $ A^{(k)}B^{(k)} 
         =  L_{n,k}   A^{\otimes k}    B^{\otimes k}   M_{p,k}$, 
    and~\eqref{eq:mixed_prod} and    Theorem~\ref{thm:multi_compi}   
     give 
    \begin{align*}
  A^{(k)}B^{(k)}&=   L_{n,k}  \left(   (AB)^{\otimes k} \right)   M_{p,k}\\&=(AB)^{(k)} . 
    \end{align*}
\end{Remark}

%%%%%%%%%%%%%%%%%%%%%
\subsection{Algebraic expression for the $k$-additive  compound}
%%%%%%%
In this section, we derive 
a formula  for the $k$-additive compound  using 
the~$k$th Kronecker sum defined in~\eqref{eq:def_k_kron_sum}.
The starting point is the following result that is an  
analogue of Theorem~\ref{thm:akk} 
to  additive compounds. 

\begin{Theorem}\label{thm:addit}
%%%%
For any~$A \in\R^{n\times n}$, and any~$k\in\{1,\dots,n\}$, we have  
\be\label{eq:asimuk}
A^{\oplus k} M_{n,k}= M_{n,k} A^{[k]} .  
\ee
%%%%%
\end{Theorem}

\begin{IEEEproof}
%%%%%%%%%%%%%%%
Fix arbitrary~$x^1,\dots,x^k \in \R^n$,  and let~$y:=\begin{bmatrix} x^1&\dots&x^k \end{bmatrix}^{(k)}$. 
Then 
    \begin{align*}
        A^{\oplus k} M_{n,k} y  
            &= A^{\oplus k} M_{n,k} \begin{bmatrix} x^{1} & \dots & x^{k} \end{bmatrix}^{(k)} \\
           &= A^{\oplus k} \sum_{(j_1,\dots,j_k) \in S(1,\dots,k)} \sigma(j_1,\dots,j_k) (x^{j_1} \otimes \dots \otimes x^{j_k}) \\
           &= \sum_{(j_1,\dots,j_k) \in S(1,\dots,k)} \sigma(j_1,\dots,j_k) \left( (Ax^{j_1} \otimes \dots \otimes x^{j_k}) + \dots + (x^{j_1} \otimes \dots \otimes Ax^{j_k}) \right) \\
%%%%%%%%%%%%%%
           &= M_{n,k} \left( \begin{bmatrix} Ax^{1} & \dots & x^{k} \end{bmatrix}^{(k)} + \dots + \begin{bmatrix} x^{1} & \dots & Ax^{k} \end{bmatrix}^{(k)} \right)  
          ,
    \end{align*}
and using~\eqref{eq:effect_k_add}  yields~$ A^{\oplus k} M_{n,k} y  = M_{n,k} A^{[k]} y$. 
\end{IEEEproof}

The next result provides a formula for the~$k$-additive compound using the~$k$th Kronecker
sum. 

 \begin{Theorem}\label{thm:kron_addi_compound}    Let $A \in \R^{n \times n}$, and fix~$k\in\{1,\dots,n\}$. Then 
    \begin{equation}\label{eq:add_kron}
        A^{[k]} = L_{n,k}A^{\oplus k } M_{n,k}. 
    \end{equation}
%%%%
\end{Theorem}
%%%%%%

\begin{IEEEproof}
   %%%
   Multiply both sides of~\eqref{eq:asimuk} from the left by~$L_{n,k}$ and use~\eqref{eq:LM_is_Identity}.
\end{IEEEproof}

\begin{Remark}
    Since~$(I_n)^{\otimes \ell}=I_{n^\ell}$, we can also write~\eqref{eq:add_kron} as
\be\label{eq:akdiff}
A^{[k]} = L_{n,k} \left( 
            A \otimes I_{n^{k-1}} 
            + I_n \otimes   A  \otimes I_{n^{k-2}} + \dots + I_{n^{k-1}}\otimes   A\right) M_{n,k}.
\ee
\end{Remark}

%%%%%%%

%%%%%%

%%%%%
\begin{Example}
    %%%%%%%%%
For~$k=2$ Eq.~\eqref{eq:add_kron} gives  \begin{align}\label{eq:add_kron_2}
        A^{[2]} &= L_{n,2} \left( 
            A \otimes I_n+   I_n\otimes   A\right) M_{n,2}\\
%%%    
            &=L_{n,2}(A\oplus A) M_{n,2}. 
        \nonumber
    \end{align}
For~$k=n$,
Eq.~\eqref{eq:add_kron} 
gives
\begin{align*}
    A^{[n]} = L_{n,n}  \left( 
            A \otimes I_n \otimes \dots \otimes I_n + I_n \otimes   A \otimes \dots \otimes I_n + \dots + I_n \otimes \dots \otimes I_n\otimes   A\right) M_{n,n},
\end{align*}
    where every Kronecker product includes~$n$ terms, 
and 
using the expressions in Example~\ref{exa:mnln}
and~\eqref{eq:mixed_prod} yields
%%%%%%%%%%%%%%%%%%%%
\begin{align*}
    A^{[n]} &= (e^1)^T \otimes \dots \otimes (e^n)^T \left( 
            A \otimes I_n \otimes \dots \otimes I_n  +\dots + I_n \otimes \dots \otimes I_n\otimes   A\right)\\&\times \sum_{(j_1,\dots,j_n) \in S(1,\dots,n)} \sigma(j_1,\dots,j_n) (e^{j_1} \otimes \dots \otimes e^{j_n}) \\
            &= ((e^1)^T A e^1) \otimes ((e^2)^T e^2) \otimes \dots \otimes ((e^n)^T e^n)
            + \dots + ((e^1)^T e^1) \otimes \dots \otimes  ((e^{n-1})^T e^{n-1}) \otimes ((e^n)^T A e^n) \\
            &= a_{11}   + \dots + a_{nn} \\
            &= \tr(A),
\end{align*}
where we used the fact that for any $(j_1,\dots,j_n) \in S(1,\dots,n) \setminus \{(1,2,\dots,n)\}$, the Kronecker product will have elements of the form $(e^i)^T e^j$ with $i \neq j$.
%%%%%%%%%
 \end{Example}

%%%%%%%%%%%%%%%%%
\begin{Remark}
%%%%%%%%%%%%%%%%%%%
    Consider the Lyapunov-like matrix ODE
    \begin{equation}\label{eq:lyap_ode}
        \dot X (t)= AX(t) + X(t) A^T,
    \end{equation}
    with $ A \in \R^{n \times n}$ and~$X:\R_+ \to\R^{n\times n}$. By the vectorization property of the Kronecker product (see~\eqref{eq:btvec}), we may rewrite~\eqref{eq:lyap_ode} in vector form as 
    \begin{equation}\label{eq:wemuli}
        \vect(\dot X) = (A \oplus A) \vect(X).
    \end{equation}
    Suppose in addition that~$X(t)$ is skew-symmetric for all~$t$. Multiplying~\eqref{eq:wemuli} on the left by $L_{n,2}$ and using 
    Remark~\ref{rem:Ln_Mn_vec} and~\eqref{eq:add_kron_2}
    gives 
    \begin{align*}
        \vech(\dot X) &= L_{n,2}(A \oplus A) \vect(X)\\
        &= L_{n,2}(A \oplus A) M_{n,2} \vech(X)\\
        &= A^{[2]} \vech(X),
    \end{align*}
    and this recovers~\cite[Eq. (8)]{Angeli2023smallgain} (see also~\cite[Section 3.1]{Angeli2014hopf}).
\end{Remark}

  %%%%%%%%%%
\section{Applications}\label{sec:appli}
 %%%%%%%%%%%%%%%%%%%%%%%%%%%%
 %%%%%%

 The formulas for compounds using Kronecker products and sums allow   using  the well-known theory of Kronecker products to analyze matrix compounds of any order~$k$. We demonstrate  this using several applications.  

 \subsection{A new proof for the expression for the entries of~$A^{[k]}$}
 %%%%%%%%%%%%%%%%%%%%%%%%%%%%%%%%%%%%%
 Theorem~\ref{thm:kron_addi_compound} can   be used to provide a simpler proof of  a known explicit expression for the entries of the additive compound~$A^{[k]}$ in terms of the entries of~$A$~\cite{schwarz1970}. To do so,   first recall  that,  by definition, every entry of the multiplicative compound~$A^{(k)}$ corresponds so some~$k\times k$ minor of~$A$, that is, every entry of~$A^{(k)}$ corresponds to rows of~$A$ with  indices  
  $(i_1,\dots,i_k) \in Q(n,k)$ and columns of~$A$ with indices~$(j_1,\dots,j_k) \in Q(m,k)$. Since~$A^{[k]}$  is defined via the $k$-multiplicative compound, the same holds for~$A^{[k]}$ as well.
%%%%%%%%%%%%%%%%
\begin{Lemma}
%%%%%%%%%%%%%%%%
    Let $A \in \R^{n \times n}$, and fix~$k \in \{1,\dots,n\}$. Let~$z$ denote the  entry of~$A^{[k]}$ corresponding to~$\alpha=(i_1,\dots,i_k) \in Q(n,k)$ and~$\beta=(j_1,\dots,j_k) \in Q(n,k)$, that is,
    $z=(A^{[k]})_{\ind^Q_{n,k}(i_1,\dots,i_k), \ind^Q_{n,k}(j_1,\dots,j_k)}$. Then 
    \begin{equation}\label{eq:zcases}
       z = \begin{cases}
            a_{i_1 i_1} + \dots + a_{i_k i_k}, & \text{if } i_p=j_p \text{ for all } p,  \\
            (-1)^{\ell+m} a_{i_\ell j_m}, & \text{if exactly  $k-1$ of the indices of  $\alpha$ coincide with those of~$\beta$, with     $i_\ell\not=j_m$ },\\
            0 , & \text{otherwise}.
        \end{cases}
    \end{equation}
\end{Lemma}
\begin{IEEEproof}
  This result was proven in~\cite{schwarz1970}. We show that  the representation  of compounds using  Kronecker products and sums  
  can be used to derive a simpler proof. 
  By Theorem~\ref{thm:kron_addi_compound},
%%%%%%%%%%
    \begin{align*}
        z &= \left(\begin{bmatrix}e^{i_1} & \dots& e^{i_k} \end{bmatrix}^{(k)}\right)^T L_{n,k} A^{\oplus k} M_{n,k} \begin{bmatrix}e^{j_1} & \dots &e^{j_k} \end{bmatrix}^{(k)} ,
    \end{align*}
    and applying~\eqref{eq:mn_alternative} gives
    %%%%%%%
    \begin{align}\label{eq:zrtp}
        z  
%%%%%%%%%
        &= (e^{i_1} \otimes \dots \otimes e^{i_k})^T (A \otimes I_n \otimes \dots \otimes I_n + \dots +I_n \otimes \dots \otimes I_n \otimes A) \nonumber  \\
        &\times \sum_{(p_1,\dots,p_k) \in S(j_1,\dots,j_k)} \sigma(p_1,\dots,p_k) (e^{p_1} \otimes \dots \otimes e^{p_k}) \nonumber \\
        &= \sum_{(p_1,\dots,p_k) \in S(j_1,\dots,j_k)} \sigma(p_1,\dots,p_k) \left( a_{i_1 p_1} \delta_{i_2 p_2} \dots \delta_{i_k p_k} + \dots + \delta_{i_1 p_1} \dots \delta_{i_{k-1} p_{k-1}} a_{i_k p_k} \right),
    \end{align}
%%%%%%%%%%%
%%%%%%%%%%%    
    where we used the mixed product property and the fact that $(e^i)^T A e^j = a_{ij}$. Note that every product of~$\delta_{ij}$s in~\eqref{eq:zrtp} includes~$k-1$ terms. 
   We consider three cases. 
   
   If~$i_p=j_p$ for~$p=1,\dots,k$ then~\eqref{eq:zrtp} gives 
   \begin{align*} 
        z  
        &= \sum_{(p_1,\dots,p_k) \in S(i_1,\dots,i_k)} \sigma(p_1,\dots,p_k) \left( a_{i_1 p_1} \delta_{i_2 p_2} \dots \delta_{i_k p_k} + \dots + \delta_{i_1 p_1} \dots \delta_{i_{k-1}  {p}_{k-1}} a_{i_k p_k} \right)\\
        &=  a_{i_1 i_1}   + \dots +  a_{i_k i_k}  . 
    \end{align*}

     If exactly  $k-1$ of the indices of  $\alpha$ coincide with those of~$\beta$, with     $i_\ell \not=j_m$\updt{,}   then the only term that may be nonzero in the sum in~\eqref{eq:zrtp} is
     \[
     \delta_{i_1p_1} \dots \delta _{i_{\ell-1}p_{\ell-1} }  a_{i_\ell p_\ell  }
       \delta_{i_{\ell+1}p_{\ell+1}} \dots \delta _{i_{k}p_{k} } , 
     \]
and this term will indeed be nonzero when~$(p_1,\dots,p_k)$
is a permutation of~$\beta=(j_1,\dots,j_k)$ such that~$p_\ell=j_m$. Since the entries of~$\beta$ and~$\alpha=(i_1,\dots,i_k)$ coincide except for~$i_\ell\not=j_m$, we conclude that in this case 
\[
z=   \sigma(p_1,\dots,p_k)   a_{i_\ell j_m }=(-1)^{\ell+m} a_{i_\ell j_m }  . 
\]
     
  If~$\alpha$ and~$\beta $ differ by two or more indices  then the fact
  that every product of~$\delta_{ij}$s in~\eqref{eq:zrtp} includes~$k-1$ terms implies that~$z=0$. This completes the proof of~\eqref{eq:zcases}. 
%%%%%%%%%%%%%%%
%%%%%%%%%%%%%%%
\end{IEEEproof}
 
Note that the proof above demonstrates how the representation of the compounds using Kronecker products and sums allows proving results about $k$-compounds, for any order~$k$,  using only properties of the Kronecker product and sum.

\subsection{Generalizing an identity   for $2\times 2$ matrices to the~$n\times n$ case}
%%%%%%%%%%%%%%%%%%%%%%%%%%
It is well-known and straightforward  to prove that for any~$A,B\in\R^{2\times 2} $, we have the identity \[
 \det(A+B)=\det(A)+\det(B)+\tr(A)\tr(B) -\tr(AB).
\]
The next result uses Theorem~\ref{thm:multi_compi} to  generalize
this identity to matrices with an arbitrary dimension~$n$. 
%%%%%%%%%%%%%%%%%%%%%%%%%%%%
\begin{Theorem}\label{thm:det_sum}
    Let~$A,B\in\R^{n\times n}$ with $n \ge 2$. Then
    \be\label{eq:aplusb_2}
    (A+B)^{(2)} = A^{(2)} + B^{(2)} + A^{[2]}B^{[2]} - (AB)^{[2]}.
    \ee
\end{Theorem}
%%%%%%%%%%%%%%%%%%%%%
\begin{IEEEproof}
%%%%%%%%%%%%%%%%
  By Theorem~\ref{thm:multi_compi},
    \begin{align*}
        (A+B)^{(2)} &= L_{n,2}(A+B)^{\otimes 2}  M_{n,2} \\
            &= L_{n,2} (A \otimes A + A \otimes B + B \otimes A + B \otimes B) M_{n,2} \\
            &= A^{(2)}   + L_{n,2}(A \otimes B + B \otimes A )M_{n,2} + B^{(2)} .
    \end{align*}
%%%%%
Using the easy  to verify  identity   
    \begin{align*} 
        A \otimes B + B \otimes A &= (A \otimes I_n + I_n \otimes A)(B \otimes I_n + I_n \otimes B) - \left( (AB) \otimes I_n + I_n \otimes (AB) \right)\\
        &= A^{\oplus 2}B^{\oplus 2} -(AB)^{\oplus 2}, 
    \end{align*}
    gives
\begin{align} 
        (A+B)^{(2)}  
            &= A^{(2)} + B^{(2)} + L_{n,2}  A^{\oplus 2}B^{\oplus 2}   M_{n,2} - L_{n,2} (AB)^{\oplus 2}  M_{n,2}    .
    \end{align}
%%%%
 Applying  Theorem~\ref{thm:addit} and Lemma~\ref{lemma:LM_IS_I} yields 
\begin{align*} 
        (A+B)^{(2)}  
            &= A^{(2)} + B^{(2)}    
             + L_{n,2}  A^{\oplus 2}M_{n,2}  B^{[2] }    - L_{n,2}  M_{n,2} (AB)^{[2]}\\
             &= A^{(2)} + B^{(2)}    
             + L_{n,2} M_{n,2} A^{[2]}   B^{[2]}    -   (AB)^{[2]}\\
             &= A^{(2)} + B^{(2)}    
             +   A^{[2]}   B^{[2] }    -   (AB)^{[2]},
    \end{align*}
 and this completes the proof. 
 %%%%
 \end{IEEEproof}

\begin{Example}
    Let $A \in \R^{n \times n}, r := \binom{n}{2}$. By Theorem~\ref{thm:det_sum},
    \begin{align*}
        (I_n+\varepsilon A)^{(2)} &=  ( I_n)^{(2)} + (\varepsilon A)^{(2)} + (I_n)^{[2]}(\varepsilon A)^{[2]} - (\varepsilon A)^{[2]}\\
        &= I_r+\varepsilon^2 A ^{(2)} +2 I_r \varepsilon A^{[2]} -\varepsilon  A^{[2]}  , 
    \end{align*}
 so   
    \begin{align*}
        \frac{d}{d\varepsilon}(I_n + \varepsilon A)^{(2)}|_{\varepsilon=0} &= \frac{d}{d\varepsilon}(I_r + \varepsilon^2 A^{(2)} +  \varepsilon   A^{[2]} )|_{\varepsilon=0} \\
            &= A^{[2]},
    \end{align*}
    as expected.
\end{Example}

%%%%%%%%%%%%%%%%%%%%%%%%%%%
\subsection{A new formula for  the $k$-additive compound of a product of matrices}
%%%%%%%%%%%%%%%%%%%%%%%%
The Cauchy-Binet Theorem asserts that~$(AB)^{(k)} = A^{(k)} B^{(k)}$, that is,   the multiplicative compound of a product of two matrices decomposes  into a product of two matrices, each depending only on one of the two original matrices. A more difficult problem is to decompose~$(AB)^{[k]}$  using terms that depend only on~$A$ or only on~$B$.
For the case~$k=2$ such a decomposition appears  in~\cite{Angeli2014hopf}, but the method there seems  difficult to extend to~$k>2$. We use the 
representation of the compounds using Kronecker products   to decompose~$(AB)^{[k]}$ in the desired form. 

 We introduce more notation. For~$n>1$, $k,i\in\{1,\dots,n\}$, and a vector~$v\in\R^n$, define a matrix~$H_{n,k,i}(v)\in\R^{(n^{k-1})\times (n^k)}$ by 
\be\label{eq:defHnki}
H_{n,k,i}(v):= I_n\otimes \dots \otimes I_n \otimes v^T \otimes I_n\otimes\dots\otimes I_n ,
\ee
where the product includes~$k$ terms,
with~$v^T$ appearing in the~$i$th 
place. For example, 
\[
H_{n,k,1}(v):= v^T \otimes  I_n  \otimes I_n\otimes\dots\otimes I_n ,
\]

\begin{Theorem}\label{thm:decomp}
%%%%%%%%%%%%%%%
Let~$A\in\R^{n\times m}$ and~$B\in\R^{m\times n}$, and fix~$k\in\{1,\dots,n \}$.
 Let~$w^i  $, $i=1,\dots,m$,
 denote the~$i$th column of~$A$, and let
    $(v^j)^T$, $j=1,\dots,m$, denote the~$j$th row of~$B$. Note that~$w^i,v^j\in\R^n$. Then
 \begin{equation}\label{eq:decomp}
     (AB) ^{[k]}   = L_{n,k}  \sum_{i=1}^m  \sum_{\ell=1}^k  (-1)^{\ell+1} H_{n,k,\ell}((w^i)^T)  
    H_{n,k,1}(v^i)    M_{n,k}.
 \end{equation}
%%%%%%%%%%%
%%%%%%%    
\end{Theorem}
%%%%%%%%%%%%%%

Proving this requires several auxiliary results. The first result shows that for any~$i $ the product~$ H_{n,k,i}(v)M_{n,k}$ is equal to~$H_{n,k,1}(v)M_{n,k}$, up to a sign change. 

\begin{Lemma}\label{lemma:H_nk}
%%%%%%%%%%%%%%%%%%%%%%%%%%%%
The matrix~$ H_{n,k,i}(v)$ satisfies 
\be\label{eq:H_MULT_M}
   H_{n,k,i}(v)M_{n,k}=(-1)^{i+1}  
H_{n,k,1}(v)M_{n,k},
\ee
%%%%%%%%%%%%%%
%%%%%%%%%%%%%%
\end{Lemma}

\begin{IEEEproof}
%%%%%%%%%5
Fix an index~$s\in\{1,\dots,k\}$. 
To simplify the notation,
denote~$Z:=H_{n,k,s}(v)M_{n,k}$.
Using~\eqref{eq:mn_alternative} and the mixed product property gives
%%%%%%%%%%%%%%%%%%%%%%
\begin{align*}
Z&=\sum_{(i_1,\dots,i_k) \in Q(n,k)} \sum_{(j_1,\dots,j_k) \in S(i_1,\dots,i_k)} \sigma(j_1,\dots,j_k) (e^{j_1} \otimes \dots \otimes e^{j_{s-1}} \otimes v^T e^{j_s} \otimes e^{j_{s+1}}\otimes \dots\otimes     e^{j_k}) \\&\times \left(\begin{bmatrix} e^{i_1} & \dots & e^{i_k} \end{bmatrix}^{(k)}\right)^T.
%%%
\end{align*}
Since~$v^T e^{j_s} = v_{j_s}$ is a scalar, we get
\begin{align*}
Z&=\sum_{(i_1,\dots,i_k) \in Q(n,k)} \sum_{(j_1,\dots,j_k) \in S(i_1,\dots,i_k)} \sigma(j_1,\dots,j_k) (  v^T e^{j_s}\otimes  e^{j_1} \otimes \dots \otimes e^{j_{s-1}}    \otimes e^{j_{s+1}}\otimes \dots\otimes     e^{j_k}) \\&\times \left(\begin{bmatrix} e^{i_1} & \dots & e^{i_k} \end{bmatrix}^{(k)}\right)^T\\
&= (v^T \otimes I_n\otimes\dots \otimes I_n)\\&\times \sum_{(i_1,\dots,i_k) \in Q(n,k)} \sum_{(j_1,\dots,j_k) \in S(i_1,\dots,i_k)} \sigma(j_1,\dots,j_k) (    e^{j_s}\otimes  e^{j_1} \otimes \dots \otimes e^{j_{s-1}}    \otimes e^{j_{s+1}}\otimes \dots\otimes     e^{j_k}) \\&\times \left(\begin{bmatrix} e^{i_1} & \dots & e^{i_k} \end{bmatrix}^{(k)}\right)^T,
%%%
\end{align*}
    and using the definition of~$M_{n,k}$ proves~\eqref{eq:H_MULT_M}.
%%%%%%%%%
%%%%%%%%%%
\end{IEEEproof}

The next result considers~$(AB)^{[k]}$ in the particular case where~$A=v$ and~$B=w^T$, with~$v,w\in\R^n$. 

\begin{Lemma}\label{lem:vw_k_add}
Let~$v,w\in\R^n$, and fix~$k\in\{1,\dots,n\}$. Then
\[
(v w ^T) ^{[k]}   =  L_{n,k}\sum_{\ell=1}^k  (-1)^{\ell+1} H_{n,k,\ell}(w^T)  
    H_{n,k,1}(v)    M_{n,k}.
\]
\end{Lemma}
\begin{IEEEproof}
    %%%%%%%%%%%%
Theorem~\ref{thm:kron_addi_compound}
gives
\begin{align*}
 (wv^T)^{[k]} 
    &=  L_{n,k} (wv^T)^{\oplus k } M_{n,k}\\
    &=  L_{n,k}(wv^T \otimes I_n \otimes \dots \otimes I_n+\dots + I_n \otimes \dots \otimes I_n\otimes wv^T ) M_{n,k},
\end{align*}
where every Kronecker product in the brackets includes~$k$ terms. Using the mixed product property gives
%%%%%%%%%%5
\begin{align}\label{eq:temp}
   (wv^T)^{[k]} 
    &=  L_{n,k}\sum_{\ell=1}^k H_{n,k,\ell}(w^T)   H_{n,k,\ell}(v)   M_{n,k}, 
\end{align}
%%%%%%%%%%%%%%
and applying Lemma~\ref{lemma:H_nk}  completes the proof. 
 %%%%%%%%%%%%%%%%%
\end{IEEEproof}

We can now prove Theorem~\ref{thm:decomp}.
%%%%%%%%%%%%%%%%%%%%%%%
\begin{IEEEproof}[Proof of Theorem~\ref{thm:decomp}]
%%%%%%%%%%%%%%%%%%%%%%%%%%%
    Let~$w^i  $ denote the~$i$th column of~$A$, and let
    $(v^i)^T$ denote the~$i$th row of~$B$. Note that~$w^i,v^i\in\R^n$. Then
$AB  =  \sum_{i=1}^m w^i (v^i)^T$, so~$(AB)^{[k]} = \sum_{i=1}^m (w^i (v^i)^T) ^{[k]}$. 
Applying Lemma~\ref{lem:vw_k_add} completes the proof. 
%%%%%%%%%%%%%%%%%
%%%%%%%%%%%%%%%%%
\end{IEEEproof}
 
The next result specializes Theorem~\ref{thm:decomp} to the case~$k=2$.
\begin{Corollary}
%%%%%%%%%%%%%%%
Let~$A\in\R^{n\times m}$ and~$B\in\R^{m\times n}$.
 Let~$w^i  $ denote the~$i$th column of~$A$, and let
    $(v^i)^T$ denote the~$i$th row of~$B$.   Then
 \[
 (AB) ^{[2]}   =    M_{n,2}^T  
 \sum_{i=1}^m          (  H_{n,2,1}(w^i) )^T  H_{n,2,1}(v^i)    M_{n,2} . 
\]
%%%%%%%%%%%
%%%%%%%    
\end{Corollary}

\begin{IEEEproof}
%%%%%%%%%%%%%%%%%%%%%%%%%%%%%%%%%%%%%%
By Theorem~\ref{thm:decomp}, 
\begin{align}\label{eq:AB2}
    (AB) ^{[2]}   &= L_{n,2}  \sum_{i=1}^m  \sum_{\ell=1}^2  (-1)^{\ell+1} H_{n,2,\ell}((w^i)^T)  
    H_{n,2,1}(v^i)    M_{n,2} \nonumber\\ 
    &=   L_{n,2}  \sum_{i=1}^m    (H_{n,2,1}((w^i)^T)   - (H_{n,2,2 }((w^i)^T) )
    H_{n,2,1}(v^i)    M_{n,2} . 
\end{align}
%%%
Let~$P:= L_{n,2}     (H_{n,2,1}((w^i)^T)   - (H_{n,2,2 }((w^i)^T) ) $. Then 
 %%% 
\begin{align*}
    P  &= 
  \sum_{(i_1, i_2) \in Q(n,2)} \begin{bmatrix} e^{i_1} &   e^{i_2} \end{bmatrix}^{(2)} (e^{i_1}  \otimes e^{i_2})^T    (  w^i  \otimes I_n - I_n \otimes  w^i   ) \\
  &=   \sum_{(i_1, i_2) \in Q(n,2)} \begin{bmatrix} e^{i_1} &   e^{i_2} \end{bmatrix}^{(2)}
     ( (e^{i_1})^T  w^i  \otimes (e^{i_2})^T  
  - (e^{i_1})^T \otimes  (e^{i_2})^T  w^i   ) \\
  & =   \sum_{(i_1, i_2) \in Q(n,2)} \begin{bmatrix} e^{i_1} &   e^{i_2} \end{bmatrix}^{(2)} 
     ( (e^{i_1})^T   \otimes (e^{i_2})^T  
  - (e^{i_2})^T \otimes  (e^{i_1})^T     )  ( w^i  \otimes I_n). 
    %%%%
\end{align*}
Thus,
\begin{align*}
P^T & = ( w^i  \otimes I_n)^T  \sum_{(i_1, i_2) \in Q(n,2)} (  e^{i_1}     \otimes  e^{i_2}    
  - e^{i_2}  \otimes   e^{i_1}     ) \left(\begin{bmatrix} e^{i_1} &   e^{i_2} \end{bmatrix}^{(2)} \right)^T\\
  &=    H_{n,2,1}(w)  M_{n,2} , 
\end{align*}
and substituting this in~\eqref{eq:AB2} completes the proof. 
 %%%%%%%%%%%%%%   
\end{IEEEproof}
%%%%%%%%%%%%%%

%%%%%%%%%%%%%%%%%%%%%%%
\section{Discussion}
%%%%%%%%%%%%%%%%%%%%%%%
Compound matrices play an important role in many fields of applied mathematics, and have recently found many applications in systems and control theory~\cite{comp_long_tutorial}. However,  the standard expressions for  these compounds are not always easy to use. For example, the compounds of a block matrix destroy the block structure. 
We derived   new formulas for matrix compounds  using    Kronecker products and   sums. This allows to use the  well-established theory of Kronecker products and sums in analyzing matrix compounds. We demonstrated several applications of these new formulas. 

An important problem 
in  systems and control theory
is deriving tractable expressions for the compounds of block matrices. 
For example, the block matrix may be the Jacobian~$J$  of a nonlinear dynamical system composed of   an  interconnection of several sub-systems,  and tractable expressions for~$J^{[k]}$ may be used to prove~$k$-contraction  based on properties of the sub-systems. 
An interesting direction  for further research is to combine our results with 
the  existing theory of  Kronecker products of block matrices 
and their  applications to large-scale systems~\cite{Hyland1989block,part_kron_prod}.

%%%%%%%%%%%%%%%%%

\section*{Appendix}
%%%%%%%%%%%%%%%%%%
We collect here several more technical results. These results   are not necessary  for deriving the main results in this paper, and are included here  for completeness.

\subsection{Expressions for~$\ind^Q$ and~$\ind^R$}
%%%%%%%%%%%%%%%%
We detail  expressions for~$\ind^Q$ and~$\ind^R$.  These expressions  are not necessary  for deriving the main results in this paper, and are included here  for completeness.
The starting point is 
  a recursive expression for~$\ind^Q_{n,k}$.
%%%%
\begin{Lemma}
    Fix an integer~$n>1$ and~$k\in\{2,\dots,n\}$. Then for any sequence~$1\leq i_1<\dots<i_k\leq n$, we have 
    \begin{equation}\label{eq:indQ_rec}
        \ind^Q_{n,k}(i_1,\dots,i_k) = \binom{n}{k} - \binom{n-i_1+1}{k} + \ind^Q_{n-i_1,k-1}(i_2 - i_1,\dots,i_k - i_1) ,
    \end{equation}
    with
    \be\label{eq:indq1}
    \ind^Q _{n,1} (i_1)= i_1.
    \ee
\end{Lemma}
Note that the solution of~\eqref{eq:indQ_rec} and~\eqref{eq:indq1} is
%%%%%%%%%%%%%%%%%%%
\begin{equation}\label{eq:indQ_nonrec}
    \ind^Q_{n,k}(i_1,\dots,i_k) = i_k + \sum_{j=1}^{k-1} \left( \binom{n - \sum_{\ell=1}^{j-1} i_\ell}{k - j + 1} - \binom{n - \sum_{\ell=1}^j i_\ell + 1}{k - j + 1} - i_j \right).
\end{equation}

\begin{IEEEproof}
%%%
Eq.~\eqref{eq:indq1} is obvious. For~$k>1$, there are in total~$\binom{n}{k}$ sequences in~$Q(n,k)$, and in~$\binom{n - i_1 + 1}{k}$ sequences  the first entry is larger than or equal to $i_1$. Hence, $\ind^Q_{n,k}(i_1,i_1+1,i_1+2,\dots,i_1+k-1) = \binom{n}{k} - \binom{n-i_1+1}{k}$. Therefore,
    \begin{equation*}
        \ind^Q_{n,k}(i_1,i_2,\dots,i_k)  = \binom{n}{k} - \binom{n-i_1+1}{k} + \ind^Q_{n-i_1,k-1}(i_2 - i_1,\dots,i_k - i_1),
    \end{equation*}
    and this completes the proof.
%%%%%%%%%%%%
\end{IEEEproof}
%%%%%%%%%%%%%%%%%%

The next result provides an    expression for~$\ind^R_{n,k}$.
   For a set~$S$, let~$|S|$ denote the cardinality of~$S$.
   
%%%%%%%%
\begin{Lemma}\label{lemma:index_R}
%%%%%%%%%%%%%%%%%%%%%%%%%%%
    Fix an integer~$n>1$ and~$k\in\{1,\dots,n\}$. Then for any sequence~$(i_1,\dots,i_k) \in R(n,k)$, we have 
    \begin{equation}\label{eq:indR}
%%%%%%%%
        \ind^R_{n,k}(i_1,\dots,i_k) = \sum_{\ell=1}^k (i_\ell - 1)n^{k-\ell} + 1.
    \end{equation}
%%%%%%
\end{Lemma}

\begin{IEEEproof}
%%%%%%%
For~$k=1$, Eq.~\eqref{eq:indR} gives~$\ind^R_{n,1}(i_1) = (i_1 - 1)n^{0} + 1 = i_1$, which is obviously   true. Suppose that~\eqref{eq:indR} is true for~$k-1$. There are~$(i_1-1)|R(n,k-1)|$ sequences in $R(n,k)$ with the first element strictly smaller than~$i_1$, so
    \begin{align*}
        \ind^R_{n,k}(i_1,\dots,i_k) &= (i_1-1)|R(n,k-1)| + \ind^R_{n,k-1}(i_2,\dots,i_k) \\
            &= (i_1 - 1)n^{k-1} + \sum_{\ell=1}^{k-1} (i_{\ell+1} - 1)n^{k-1-\ell} + 1 \\
            &= \sum_{\ell=1}^k (i_\ell - 1)n^{k - \ell} + 1,
    \end{align*}
    where we used the fact that the cardinality of $R(n,k)$ is $|R(n,k)| = n^k$.
\end{IEEEproof}

For example, for~$k=n$ and the sequence~$(1,2,\dots,n)$ Eq.~\eqref{eq:indR}
gives
%%%%%%%%
\begin{align}
        \ind^R_{n,n}(1,\dots,n)& = \sum_{\ell=1}^n ( \ell - 1)n^{n-\ell} + 1\nonumber\\
        &=\frac{n^n-n}{(n-1)^2}.
%%%%%%%%%%%%
\end{align}

\bibliographystyle{IEEEtranS}
\bibliography{refs}

\end{document}